\documentclass[10pt]{article}

\usepackage{latexsym}
\usepackage{enumerate, amssymb, amsmath, amsthm}
\usepackage[utf8]{inputenc}
\usepackage{graphicx}
\usepackage{url}

\theoremstyle{theorem}

\theoremstyle{definition}

\begin{document}

\title{Envelopes are solving machines for quadratics and cubics and certain polynomials of arbitrary degree}
\author{Michael Schmitz and André Streicher}

\maketitle

Everybody knows from school how to graphically solve a quadratic equation
$$x^2-px+q=0,$$
if $p,q \in \mathbb{R}$ are given. Simply plot the graph of $f(x)=x^2-px+q$ and find the point(s) of intersection with the $x$-axis. If $p$ and $q$ are modified, you have to start over and draw a new parabola to find the solutions. Can't there be one single curve that simultaneously solves all quadratic equations?

Stunningly, there actually is such a magic parabola that serves as a solving machine for any quadratic equation, namely $f(x) = \frac{1}{4} x^2$. Of course, the solutions are no longer given by points of intersection with the $x$-axis, but are obtained by drawing tangent lines to $f$ through a given point $(p,q)$. Moreover, the technique can be generalized to equations of the form $x^n - px + q = 0$, and the number of real solutions of such an equation can be seen immediately.

In this article, which is strongly inspired by lecture 8 from the wonderful book \cite{Fuchs_Tabachnikov}, we derive the above mentioned methods in an elementary way and conclude by pointing out relations to the duality of points and lines in the plane and the concept of Legendre transformation.

\section*{Solving quadratic equations}

We consider an ordinary quadratic equation $x^2 -px + q = 0$\footnote{We use the form $x^2-px+q=0$ instead of $x^2+px+q=0$ only for convenience reasons; due to that, most calculations in the further course look nicer.} with $p,q \in \mathbb{R}$. Of course, everybody knows from school how to solve it, but we want to look at it from an unusual perspective. We solve for $q$ obtaining $q=xp-x^2$ and interpret $q$ as a function of $p$, depending on the parameter $x$. This gives a family of linear functions
$$Q_x(p) = xp -x^2,$$
where $Q_x$ has slope $x$ and axis intercept $-x^2$. Picture 1 shows $Q_1$ and $Q_2$.
\begin{center}
	\includegraphics[width=0.5\textwidth]{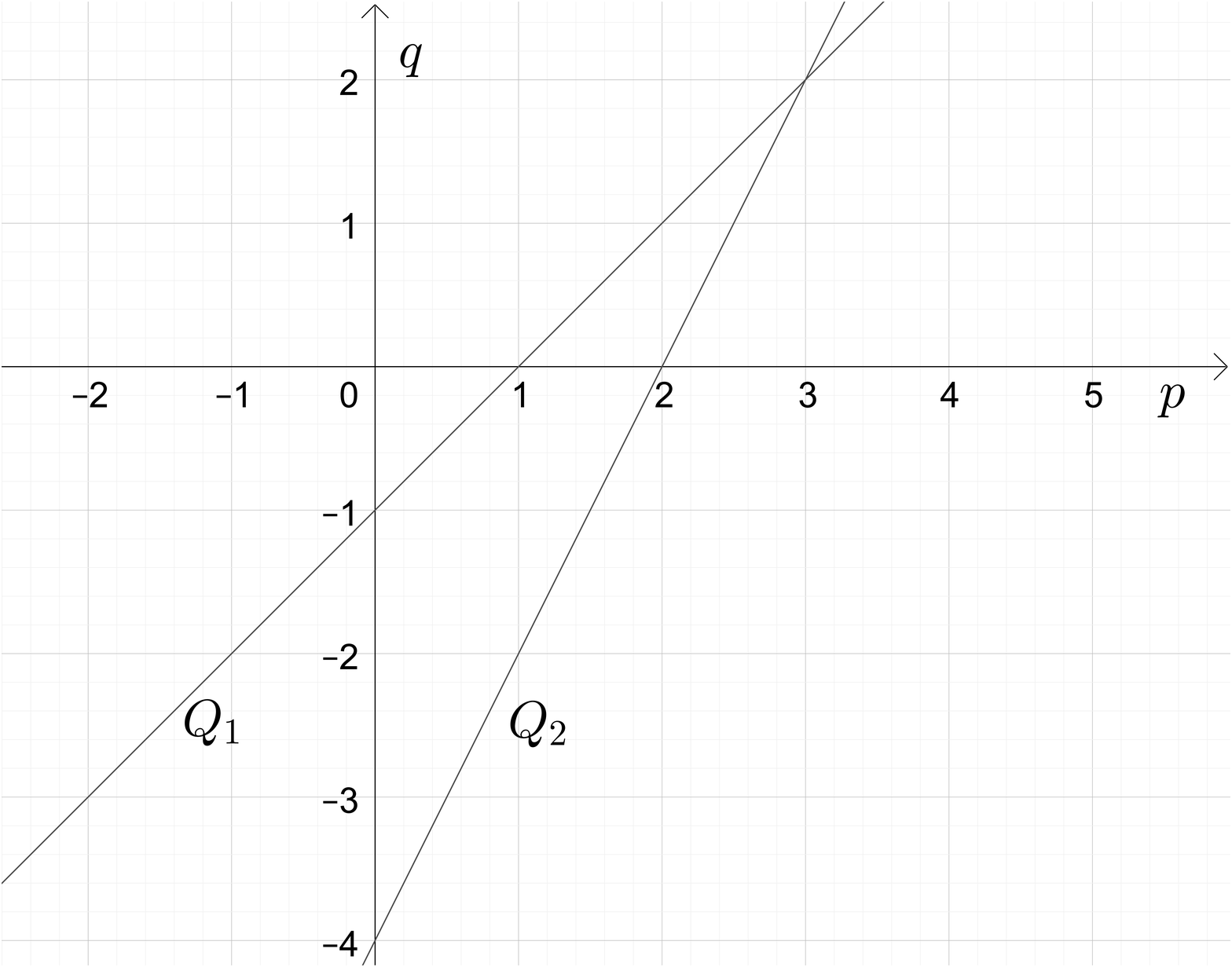}\\
	{\small Picture 1: the lines $Q_1$ and $Q_2$}
\end{center}

The graph of $Q_1$ describes all quadratic equations having $x=1$ as a solution, i.\:e., $x=1$ is a solution of $x^2-px+q=0$ if and only if $(p,q)$ lies on $Q_1$. Therefore, the unique pair $(p,q)$ of parameters of a quadratic equation with solution set $\{1,2\}$ should be given by the point of intersection of $Q_1$ and $Q_2$. In fact, the point of intersection is $(3,2)$ and the corresponding equation is $x^2 -3x+2=0$.

As a quadratic equation has at most two solutions, not more than two of the lines $Q_x$ can go through one given point $(p,q)$. Nevertheless, we want to derive this fact from the mere form of the lines $Q_x$. Therefore, we consider $x\ne y$ such that $Q_x$ and $Q_y$ go through $(p,q)$. We show that if $(p,q)$ lies also on $Q_z$, then it follows that $z=x$ or $z=y$.

The point $(p,q)$ lies on $Q_x$ if and only if $q=xp-x^2$. Thus, that $Q_x$ and $Q_y$ intersect in $(p,q)$ implies $xp-x^2 = yp - y^2$, from which we obtain $p(x-y) = x^2-y^2$. Division by $x-y$ ($\ne0$) yields $x+y = p$.\footnote{Note that we hereby obtained one part of Vieta's formula for a quadratic. We will investigate this later in more detail.} Now, if $Q_z$ goes also through $(p,q)$ and $x\ne z$, it follows analogously that $x+z = p$. Subtracting these equations gives $y=z$.\\

Considering picture 2, which shows a few more of the lines $Q_x$, leads to the conjecture that they have a quadratic envelope\footnote{By envelope we mean a differentiable function $e$ such that for every $p$ there is one unique $x$ for which $Q_x$ is a tangent to $e$ at the point $(p,e(p))$.}, and we are going to determine it.
\begin{center}
	\includegraphics[width=0.7\textwidth]{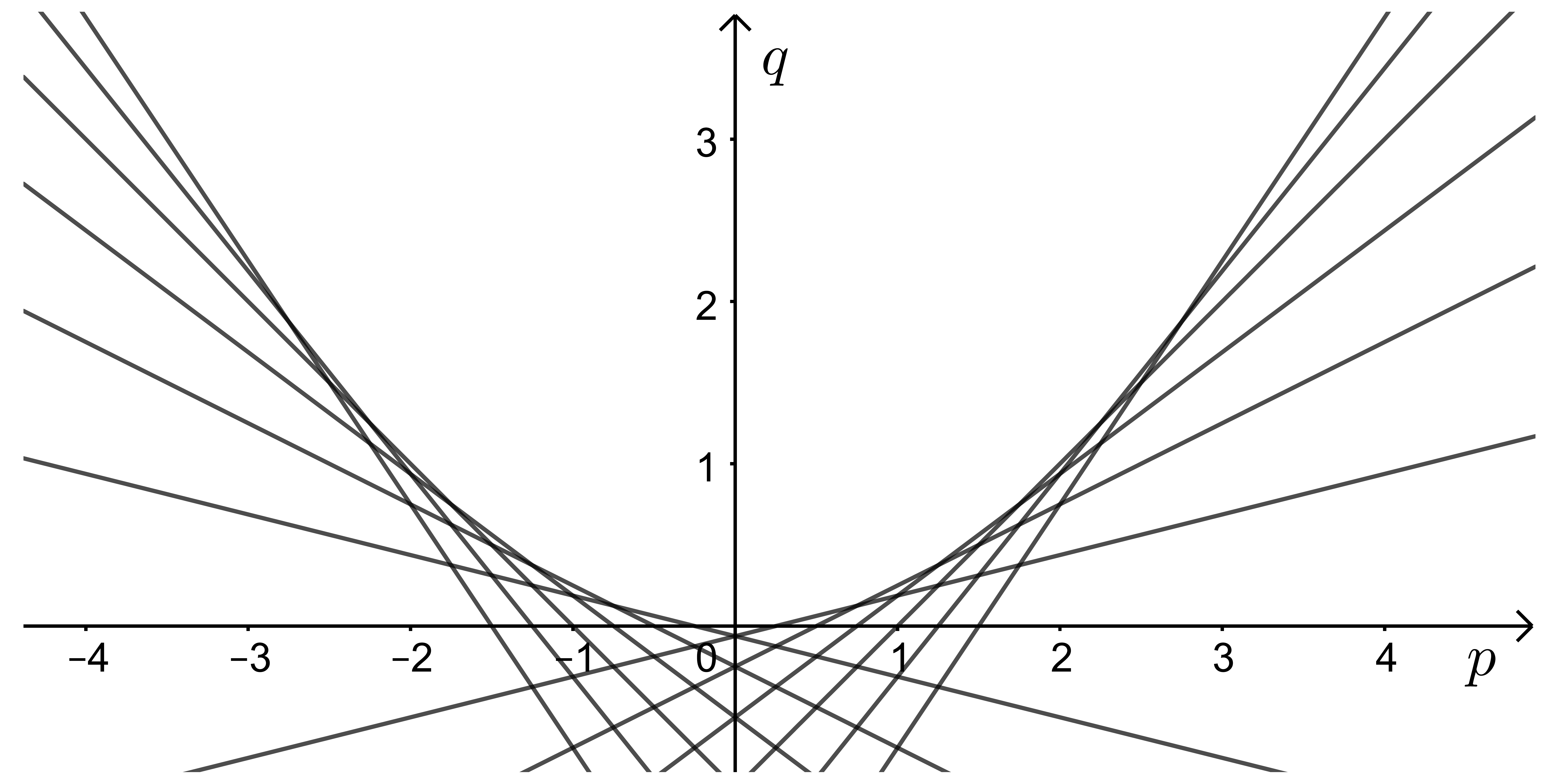}\\
	{\small Picture 2: some lines from the family $\{Q_x ~|~x \in \mathbb{R}\}$}
\end{center}

To this end, we find the point of intersection $p_\varepsilon$ of two lines $Q_x$ and $Q_{x+\varepsilon}$ for some small $\varepsilon \ne 0$ and determine the limit as $\varepsilon \to 0$. We have
$$Q_{x+\varepsilon}(p_\varepsilon) = Q_x(p_\varepsilon) \Leftrightarrow (x+\varepsilon)p_\varepsilon - (x+\varepsilon)^2 = x p_\varepsilon - x^2,$$
and the latter is equivalent to $p_\varepsilon = 2x + \varepsilon$. Thus, we obtain $p_\varepsilon \to 2x$ as $\varepsilon \to 0$, let $p=2x$ be the `point of intersection of two infinitesimal distinct lines' from our family, and denote the envelope function by $e$. It follows that $x=\frac{p}{2}$, and
$$e(p) = Q_{\frac{p}{2}}(p) = \frac{p}{2} \cdot p - \left(\frac{p}{2}\right)^2 = \frac{p^2}{4}.$$
Now that we know $e(p)$ it can be easily used to graphically solve any quadratic equation $x^2-px+q=0$. Given $p,q$ we just have to find the line(s) $Q_x$ on which the point $(p,q)$ lies. Therefore we construct the tangent line(s) to $e$ through $(p,q)$. Picture 3 illustrates the method by the example of solving $x^2-x-2=0$. For convenience we marked the rescaling $x=\frac{p}{2}$ on the graph of $e$, so that one can simply read off the solutions from the picture. 
\begin{center}
	\includegraphics[width=0.9\textwidth]{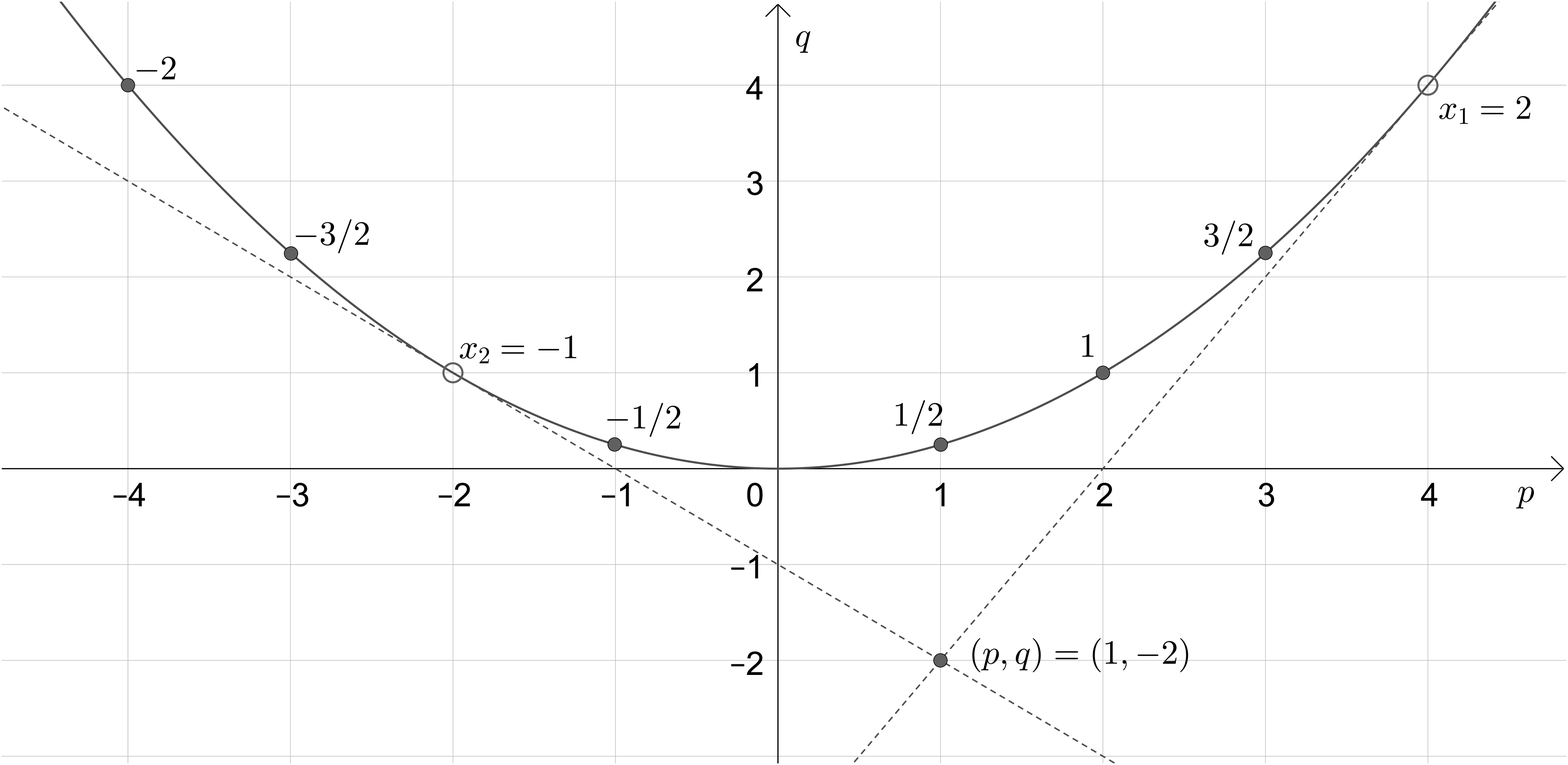}\\
	{\small Picture 3: solving $x^2-x-2=0$}
\end{center}

Isn't it wonderful? We can solve any quadratic equation with one picture. After discovering this `magic parabola' with a class, one could draw a large version of it on a whiteboard and construct the tangent lines with ropes and magnets. Of course, also a GeoGebra applet is easily constructed and fun.\\

Moreover, picture 3 tells us the number of solutions of any given quadratic equation. By convexity there are two tangent lines (i.\:e. two solutions of the equation) to the graph of $e$ through a given point $(p,q)$ if and only if $(p,q)$ lies `below' $e$. There is exactly one tangent line if and only if $(p,q)$ lies on the graph of $e$, and there is no tangent line if and only if $(p,q)$ lies `above' $e$. This can be expressed algebraically by noting that $(p,q)$ lies on [below, above] $e$ if and only if $q=e(p)$ [$q<e(p), q>e(p)$]. As
$$q=e(p) \Leftrightarrow \frac{p^2}{4}-q = 0,$$
we have rediscovered the well known discriminant of a quadratic, which indicates the number of solutions.

\section*{Solving cubic equations}

As a next step we want to extend our method to cubic equations of the form
\begin{equation}\label{cubic_equation}
	x^3-px+q=0.
\end{equation}
 In the further course we will deal with such equations for arbitrary $n$ instead of 3, and of course, the results of this and the previous section follow from the general considerations. Nevertheless, we think it is worth doing these special cases first, because they are illustrative for the case study ($n$ even or odd) that is needed later.

 Following the same idea as above we rearrange to $q=xp - x^3$ and interpret $q$ as a function of $p$ with parameter $x$. That is, we consider the family of linear functions
$$Q_x(p) = xp - x^3,$$
where $Q_x$ has slope $x$ and axis intercept $-x^3$.

Picture 4 shows some of these lines, and we get the impression that there exists an envelope with two branches, which we like to determine.

\begin{center}
	\includegraphics[width=0.7\textwidth]{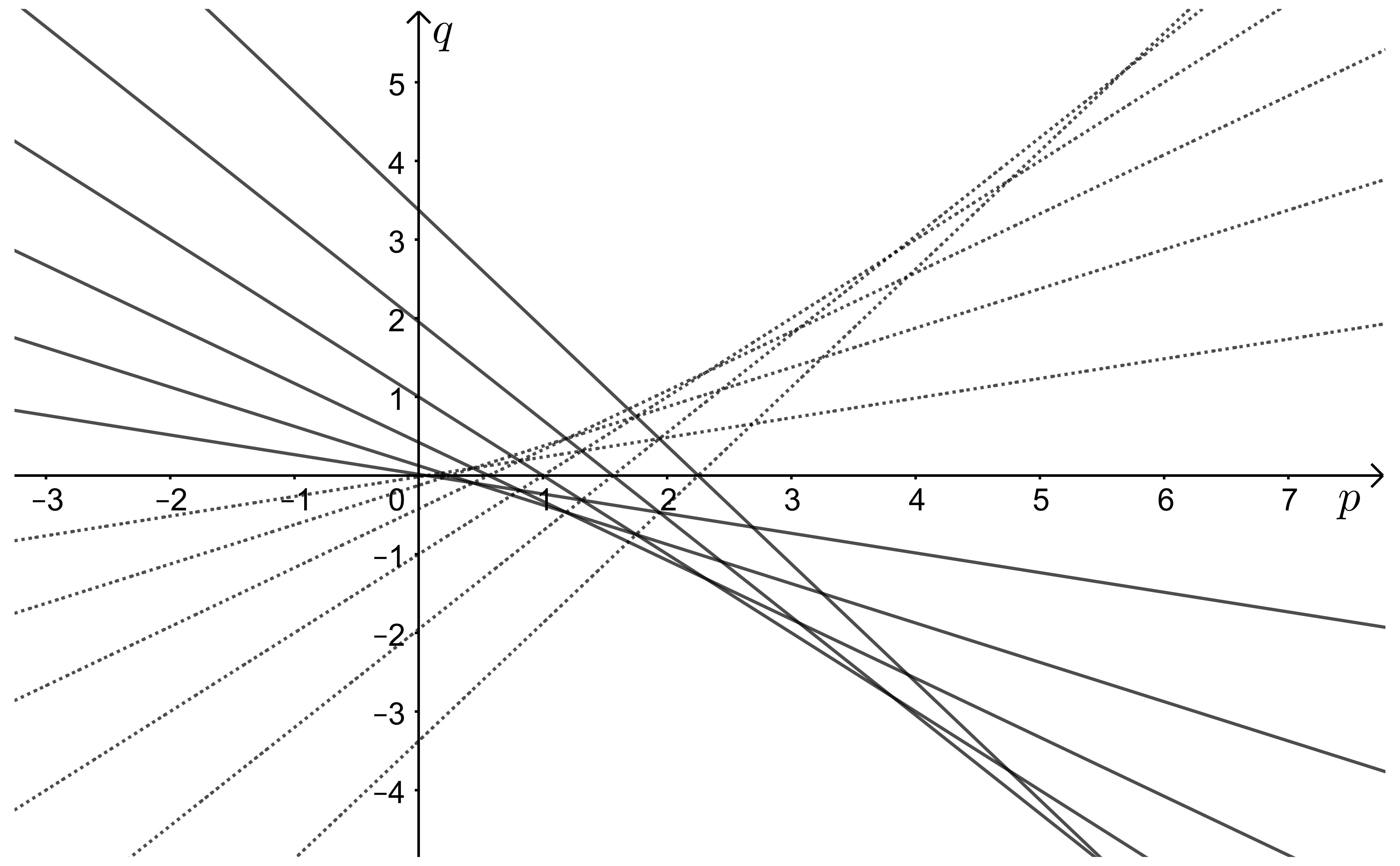}\\
	{\small Picture 4: some lines from the family $\{Q_x ~|~ x \in \mathbb{R}\}$}
\end{center}

As before, we find the point $p_{\varepsilon}$ of intersection of two lines $Q_x$ and $Q_{x+\varepsilon}$ for some small $\varepsilon \ne 0$ and let $\varepsilon \to 0$. We have
$$Q_{x+\varepsilon}(p_\varepsilon) = Q_x(p_\varepsilon) \Leftrightarrow (x+\varepsilon)p_\varepsilon - (x+\varepsilon)^3 = xp_\varepsilon - x^3.$$
The last equation is equivalent to
$$\varepsilon p_\varepsilon = (x+\varepsilon)^3 - x^3.$$
Expanding on the right hand side and dividing by $\varepsilon$ on both sides yields
$$p_\varepsilon = 3x^2 + 3x \varepsilon + \varepsilon^2,$$
and we see that $p_\varepsilon \to 3x^2$ as $\varepsilon \to 0$. Again, we interpret $p=3x^2$ as `infinitesimal point of intersection' and obtain that $p \ge 0$. Solving for $x$ gives $x = \pm \sqrt{p/3}$, and we see that we will actually get two envelope branches. For the first branch we obtain
$$e(p) = Q_{\sqrt{\frac{p}{3}}}(p) = \sqrt{\frac{p}{3}} \cdot p - \left( \sqrt{\frac{p}{3}} \right)^3 = 2 \left( \frac{p}{3} \right)^{\frac{3}{2}},$$
and the second branch equals $-e(p)$. Picture 5 shows both branches of the envelope.

\begin{center}
	\includegraphics[width=0.7\textwidth]{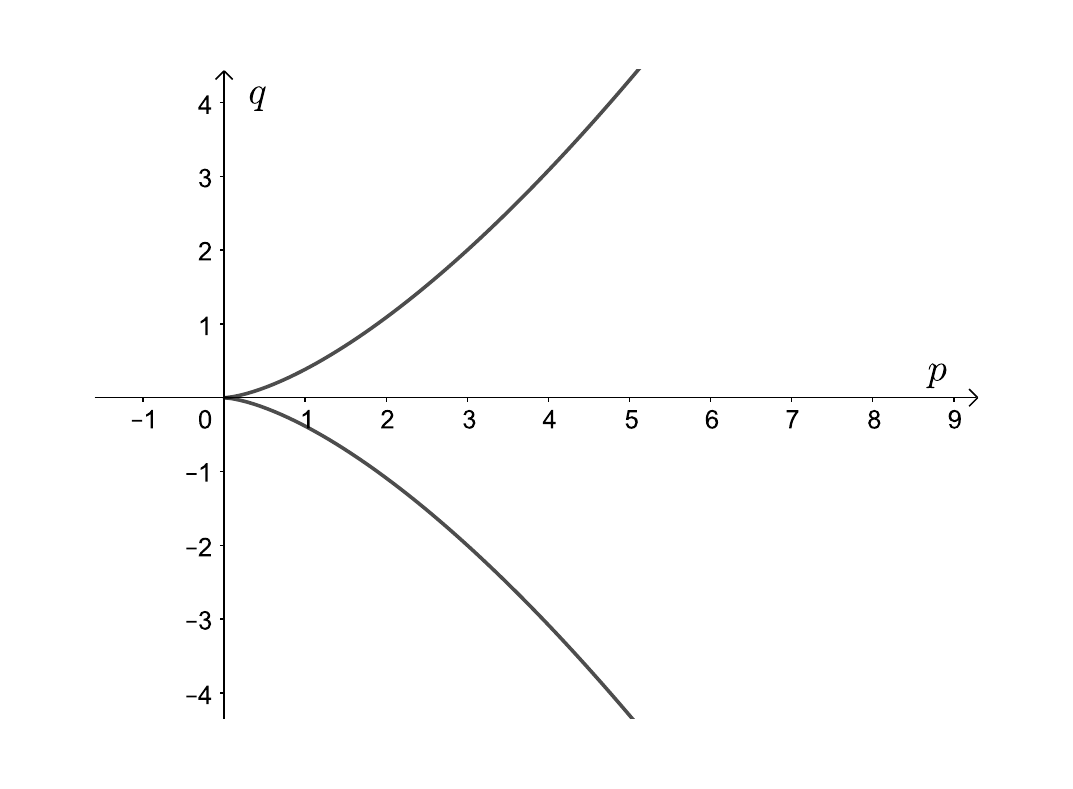}\\
	{\small Picture 5: a two-branched envelope}
\end{center}

To use the envelope for solving cubic equations we have to find all tangents to these branches through a given point $(p,q)$. Of course, the number of tangents (or solutions) depends on the position of the point.

Both branches are continuously differentiable on $\mathbb{R}_{>0}$ with $e'(p)=\sqrt{p/3}$ and differentiable from the right at the origin with $e'(0)=0$. Furthermore, we have $\lim_{p \to \infty} e'(p) = \infty$, and therefore $e$ takes on every slope in $[0,\infty)$, whereas $-e(p)$ takes on every non-positive slope. $e$ is strictly convex and $-e$ is strictly concave. From these considerations we can conclude the following.

\begin{itemize}
	\item There are exactly three different tangents through $(p,q) \ne (0,0)$ to the envelope if and only if $(p,q)$ lies in the region strictly `between the branches' (see picture 6), because in this case there are two tangents to one branch and one tangent to the other branch. A special case is when $(p,q)$ lies on the $x$-axis. Then there are two tangents to both branches, but the axis itself is an identical tangent.
	\item There are exactly two different tangents through $(p,q) \ne (0,0)$ to the envelope if and only if $(p,q)$ lies on one of the branches, because in this case there is a unique tangent to both branches.
	\item There is one unique tangent through $(p,q) \ne (0,0)$ to the envelope if and only if $(p,q)$ lies in the region strictly `not between the branches' (see picture 6), because in this case there is no tangent to one of the branches and one unique tangent to the other one. A special case is when $(p,q)$ lies on the $x$-axis. Then the axis itself is the only mutual tangent.
\end{itemize}

\begin{center}
	\includegraphics[width=0.9\textwidth]{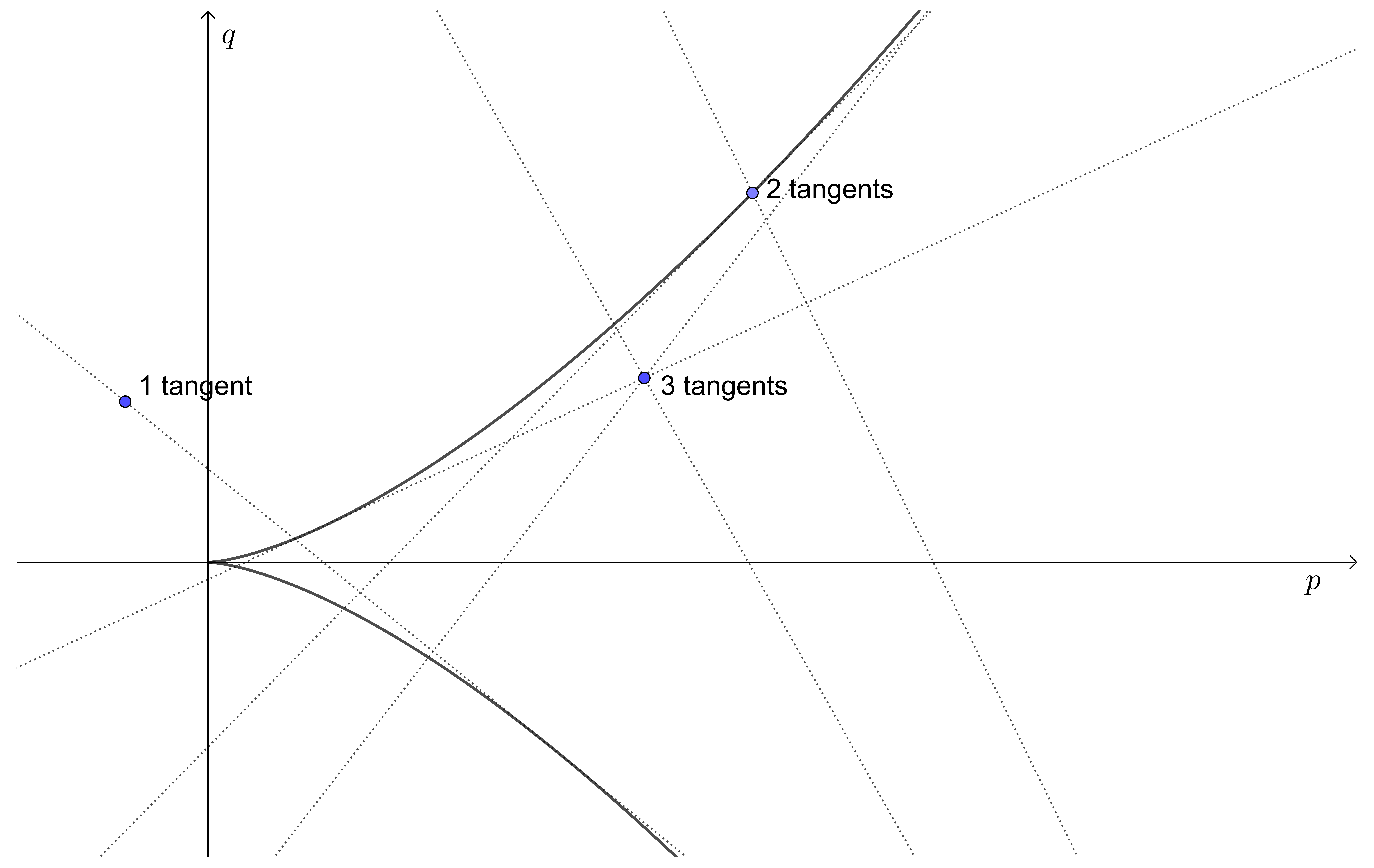}\\
	{\small Picture 6: different positions of $(p,q)$}
\end{center}

\textbf{The discriminant of a cubic equation.} The point $(p,q) \ne (0,0)$ lies on the envelope if and only if $|q| = e(p)$. It lies strictly between the branches if the analogous condition holds with $<$ instead of $=$, and strictly not between the branches if we have the inequality with $>$. Therefore, we obtain that equation (\ref{cubic_equation}) has exactly two solutions if and only if $q^2 = 4 (p/3)^3$ and this equivalent to
$$\Big(\frac{q}{2}\Big)^2 - \Big(\frac{p}{3}\Big)^3 = 0.$$
For one and three solutions we have analogous inequalities, and hereby rediscovered the well known discriminant of a cubic equation without quadratic term.\footnote{By a simple substitution of the form $x \to x+c$ every cubic equation can be transformed into one without quadratic term. For more details on this and various considerations on cubic equations from a mathematics-educational perspective see \cite{Schmitz}.}

\section*{Solving equations of arbitrary degree}

We want to generalize our ideas to equations of the form
\begin{equation}\label{general_equation}
x^n - px + q = 0,
\end{equation}

where $n\ge 2$ is an arbitrary integer. As before, the equation defines a family of linear functions with parameter $x \in \mathbb{R}$:
$$Q_x(p) = xp - x^n$$
To determine the envelope we use the same approach as above, namely
\begin{eqnarray*}
	Q_{x+\varepsilon}(p_\varepsilon) = Q_x(p_\varepsilon) &\Leftrightarrow& (x+\varepsilon)p_\varepsilon - (x+\varepsilon)^n = xp_\varepsilon - x^n\\
	&\Leftrightarrow&  \varepsilon p_\varepsilon = \sum_{k=0}^{n} \binom{n}{k} x^k \varepsilon^{n-k} - x^n\\
	&\Leftrightarrow& p_\varepsilon =  \sum_{k=0}^{n-1} \binom{n}{k} x^k \varepsilon^{n-k-1}.
\end{eqnarray*}
Letting $\varepsilon \to 0$ each term including a factor $\varepsilon$ with a positive exponent vanishes. As the only summand with a non-positive $\varepsilon$-exponent is the $(n-1)$th, it follows that $p_\varepsilon \to \binom{n}{n-1}x^{n-1} = n x^{n-1}$ for $\varepsilon \to 0$. Thus, we obtain $$p=n x^{n-1}$$
as `infinitesimal point of intersection' and distinguish two cases to proceed.\\

\textbf{Case 1: $n$ is even.} Then $n-1$ is odd, we obtain $x=\left( \frac{p}{n}\right)^{\frac{1}{n-1}}$, and therefore
$$e(p) = Q_{\left(\frac{p}{n}\right)^{\frac{1}{n-1}}}(p) = \left(\frac{p}{n}\right)^{\frac{1}{n-1}} \cdot p - \left(\frac{p}{n}\right)^{\frac{n}{n-1}} = (n-1) \left(\frac{p}{n}\right)^{\frac{n}{n-1}}.$$
The following picture shows the envelopes for $n=2,4,6$ (dotted, dashed, solid).
\begin{center}
	\includegraphics[width=0.7\textwidth]{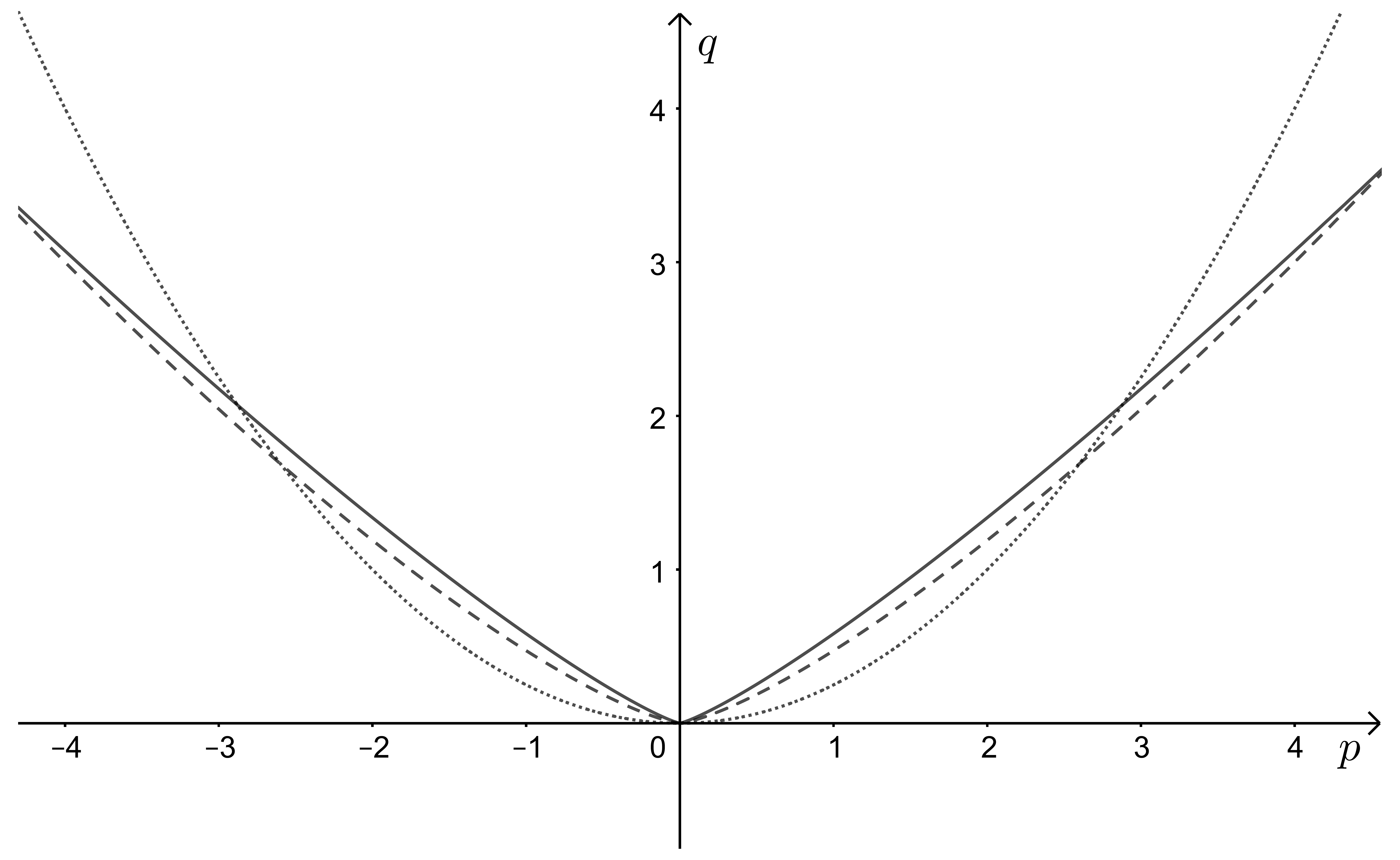}\\
	{\small Picture 7: $e(p)$ for $n=2,4,6$}
\end{center}

From the picture one might get the impression that $e$ is no longer differentiable at the origin for larger $n$, but this is not the case. $e$ is continuously differentiable on the whole domain $\mathbb{R}$ with $e'(p)=\left(\frac{p}{n}\right)^{\frac{1}{n-1}}$, in particular $e'(0)=0$. Furthermore, we see $e'(p) \to \infty$ for $p \to \pm \infty$. Therefore $e$ takes on every slope in $[0,\infty)$ and we can conclude: Through a point $(p,q)$ there exist(s)
\begin{itemize}
	\item exactly two tangents to $e$ if $(p,q)$ lies below $e$, that is, if and only if $q<e(p)$, and in this case equation (\ref{general_equation}) has exactly two solutions.
	\item one unique tangent to $e$ if $(p,q)$ lies on $e$, that is, if and only if $q=e(p)$, and in this case equation (\ref{general_equation}) has one unique solutions.
	\item no tangent to $e$ if $(p,q)$ lies above $e$, that is, if and only if $q>e(p)$, and in this case equation (\ref{general_equation}) has no solutions.
\end{itemize}

\textbf{Case 2: $n$ is odd.} Then $n-1$ is even, and since $p=nx^{n-1}$ we obtain that $p \ge 0$. We have $x = \pm (p/n)^{\frac{1}{n-1}}$, so that the envelope will have two branches. Similar to case 1 we obtain the first branch by
$$e(p) = Q_{\left(\frac{p}{n}\right)^{\frac{1}{n-1}}} (p) = \left(\frac{p}{n}\right)^{\frac{1}{n-1}} \cdot p - \left(\frac{p}{n}\right)^{\frac{n}{n-1}}
= (n-1) \left(\frac{p}{n}\right)^{\frac{n}{n-1}},$$
and the second branch equals $-e(p)$. Picture 8 shows both envelope branches for $n=3,5,7$ (dotted, dashed, solid).
\begin{center}
	\includegraphics[width=0.7\textwidth]{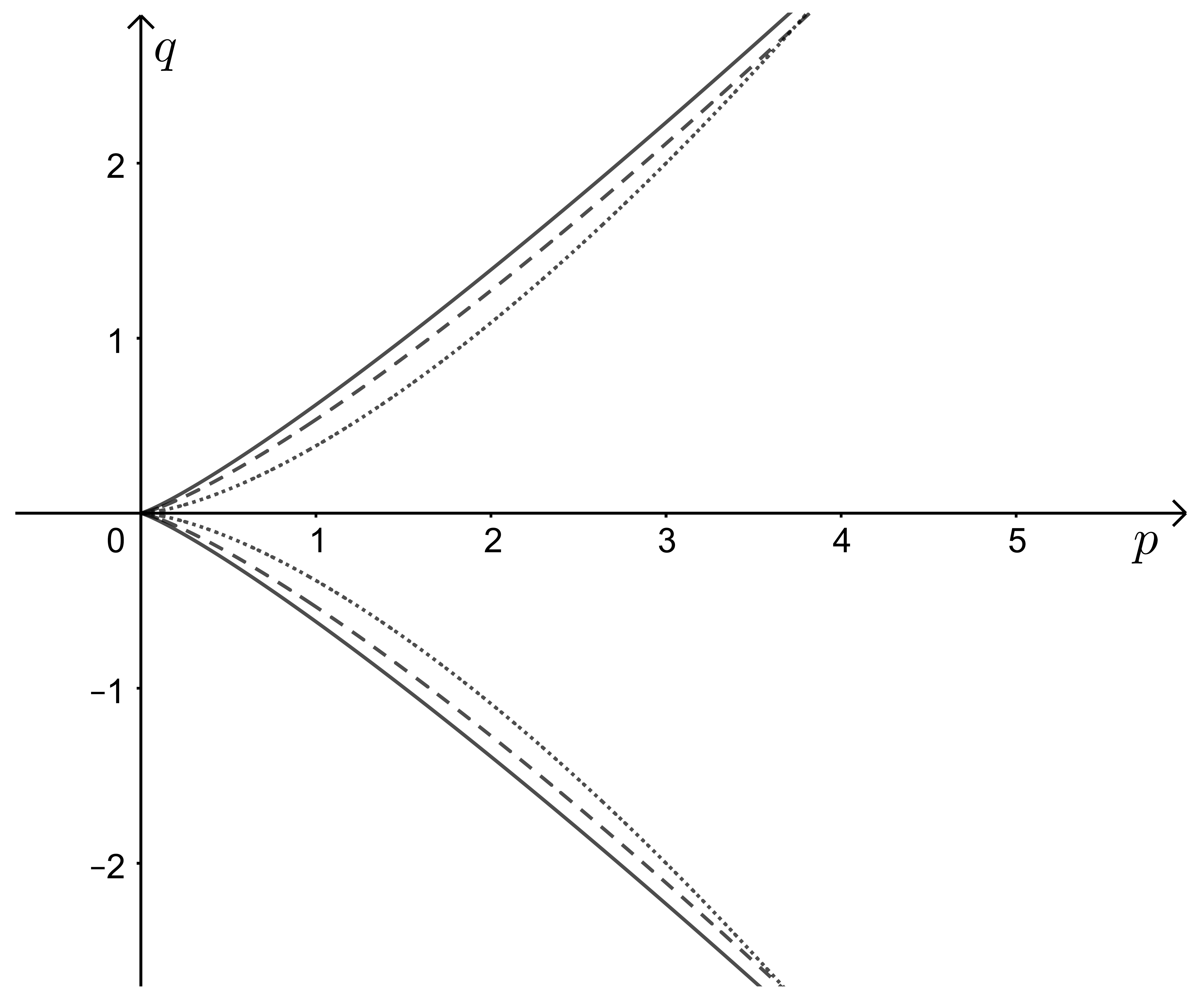}\\
	{\small Picture 8: $e(p)$ for $n=3,5,7$}
\end{center}

The branch $e$ (and of course also $-e$) is continuously differentiable on $\mathbb{R}_{>0}$, and we have $e'(p)=\left(\frac{p}{n}\right)^{\frac{1}{n-1}}$. At the origin $e$ and $-e$ are differentiable from the right and we have $e'(0)=-e'(0)=0$. Furthermore, it holds true that $e'(p)>0$ for $p>0$, $e'(p) \to \infty$ for $p \to \infty$, and $-e'(p) \to -\infty$ für $p \to \infty$. Therefore, $e$ takes on every slope in $[0,\infty)$ and $-e$ takes on every slope in $(-\infty,0]$. Thus, we obtain that through a point $(p,q) \ne (0,0)$ there exist(s)

\begin{itemize}
	\item exactly three tangents to the envelope if $(p,q)$ lies in the region `strictly between' the branches, that is, if and only if $|q|<e(p)$.
	\item exactly two tangents to the envelope if $(p,q)$ lies on one of the branches, that is if and only if $|q|=e(p)$.
	\item one unique tangent to the envelope if $(p,q)$ lies in the region strictly `not between' the branches, that is if and only if $|q|>e(p)$.
\end{itemize}

\textbf{Summary of both cases and the discriminant.} For all $n$ it holds true that $e(p)=(n-1) \left( \frac{p}{n} \right)^{\frac{n}{n-1}}$. If $n$ is even, the envelope has one branch and is defined on the entire real line. If $n$ is odd, the envelope has the two branches $e$ and $-e$, and both of them are defined on $\mathbb{R}_{\ge 0}$. For both, odd and even $n$, the relation ($<,>$ or $=$) between $q$ (or $|q|$) and $e(p)$ indicates the number of solution of equation (\ref{general_equation}). From this we can derive the determinant of (\ref{general_equation}). For even $n$ it holds true that
$$q=e(p) \Leftrightarrow \frac{q}{n-1} = \left( \frac{p}{n} \right)^{\frac{n}{n-1}} \Leftrightarrow \Big( \frac{p}{n} \Big)^n -\Big( \frac{q}{n-1} \Big)^{n-1} = 0.$$
The calculation for odd $n$ is analogous, so the determinant of equation (\ref{general_equation}) is given by $\Big( \frac{p}{n} \Big)^n -\Big( \frac{q}{n-1} \Big)^{n-1}$. Isn't it aesthetic? 

\section*{Final remarks and prospect}

\subsection*{Duality} It is worth mentioning that we used a more general concept here, namely the duality between straight lines and points in the plane. On the one hand, a linear equation of the form $q=mp+n$ represents a line in the $pq$-plane, which is uniquely determined by the pair $(m,n)$ of slope and axis intercept. On the other hand, this pair represents a point in the $mn$-plane. Thus, we have a one-to-one correspondence between the (non-vertical) straight lines in the $pq$-plane and the points in the $mn$-plane.

If we consider a linear function $Q(p) = mp + n$ with slope $m$ and axis intersect $n$, a point $(p,q)$ in the $pq$-plane lies on the graph of $Q$ if and only if $q=mp+n$, which is equivalent to $n = -pm  + q$. This means that in the $mn$-plane the point $(m,n)$, which is corresponding to $Q$, lies on the graph of  $N(m) = -pm + q$ with slope $-p$ and axis intercept $q$.

To illustrate this, picture 9 shows the lines $a$, given by $q=-p+3$ and $b$, given by $q=p-1$, which intersect in the point $S=(2,1)$. The corresponding points $A=(-1,3)$ and $B=(1,-1)$ in the $mn$-plane lie on the line $s$, given by $n=-2m+1$.

\begin{center}
	\includegraphics[width=0.425\textwidth]{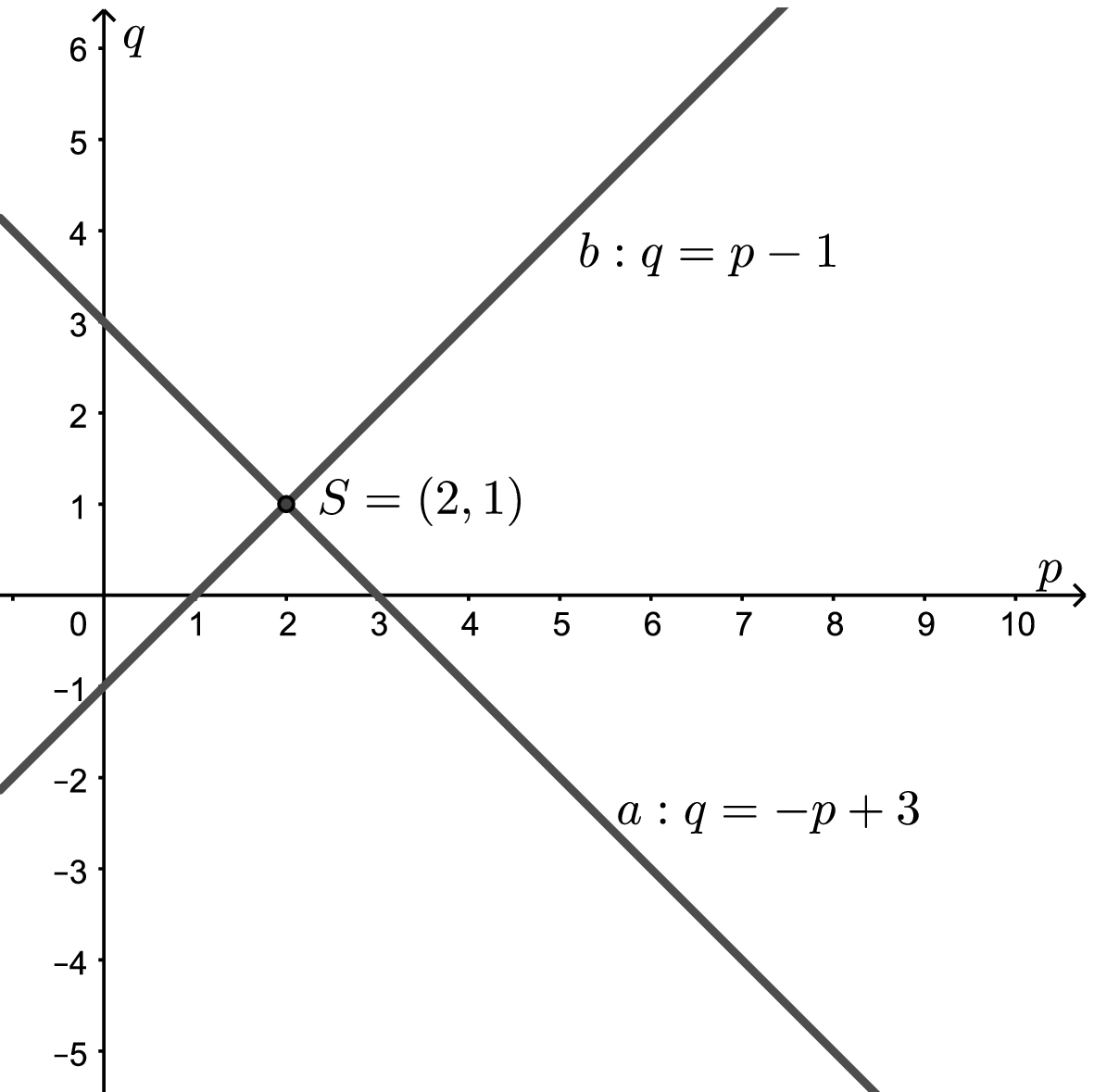}\hspace{6mm}
	\includegraphics[width=0.452\textwidth]{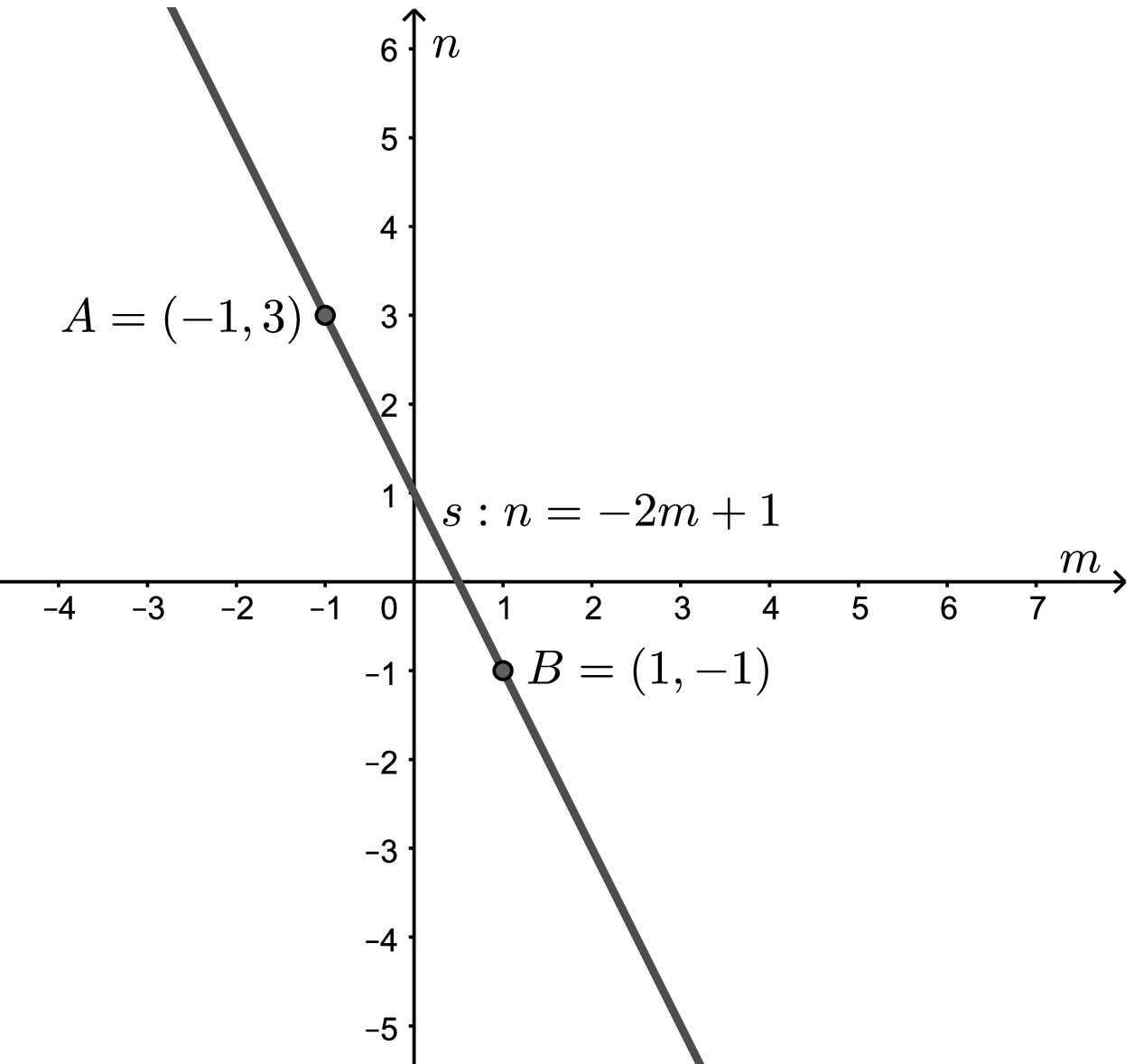}\\
	{\small Picture 9: duality between lines in the $pq$-plane and points in the $mn$-plane}
\end{center}

Duality can sometimes be helpful, which we want to illustrate by deriving Vieta's formula for quadratic equations in a nice and uncommon way. We assume that $x^2-px+q = 0$ has the solutions $u$ and $v$. Then $(p,q)$ is the point of intersection of the lines $Q_u(p) = up - u^2$ and $Q_v(p) = v p - v^2$. By duality the corresponding points $(u,-u^2)$ and $(v,-v^2)$ lie on the line with slope $-p$ and axis intercept $q$.

As the line through these two points has slope $\frac{u^2-v^2}{v-u} = -(u+v)$, we obtain $p=u+v$. Moreover, by plugging in $(u,-u^2)$, the axis intercept $q$ is obtained:
$$q = pu - u^2 = (u+v)u - u^2 = uv.$$

\subsection*{Relation to Legendre transforms} Readers that are familiar with Legendre transforms might have noticed that our magic envelope matches with the Legendre transform of $f(x)=x^n$ (at least for even $n$, so that $f$ is convex). Concluding this article we want to point out that this is not a coincidence.

If we consider any smooth and strictly convex function $f:I \to \mathbb{R}$, where $I$ is an interval of reals, we can represent $f$ in terms of its first derivative in the following way: Because $f$ is strictly convex, we have $f''(x)>0$ for all $x \in I$. Therefore, $f'$ is strictly increasing on $I$ and hence a one-to-one-function. Thus, for every $p \in I^*:=\{f'(x)~|~ x \in I\}$ there is exactly one $x$ with $p=f'(x)$, or -- geometrically expressed -- no two tangent lines to $f$ have the same slope.

Now, instead of $I$ we can use the set $I^*$ of slopes $p$ as the domain of a new function that contains all information about $f$. The tangent line $t$ to $f$ at a given point $(x_0,f(x_0))$ has slope $p:=f'(x_0)$ and determines its $y$-axis intercept uniquely. The Legendre transform $f^*$ of $f$ maps $p$ to the negative\footnote{It is a convention to take the negative.} of the axis intercept of $t$. A formula for $f^*(p)$ is easily derived. We have
$$t(x) = f'(x_0)(x-x_0) + f(x_0) = px - (px_0 - f(x_0)),$$
and therefore the axis intercept of $t$ equals $px_0 - f(x_0)$. For a given function $f$ we can express this solely in terms of $p$ using the one-to-one relation between $p=f'(x_0)$ and $x_0$.

For instance, let us consider $f(x)=x^2$. Let $p=f'(x) = 2x$ be the slope of a tangent line to $f$ at some point $(x,f(x))$.\footnote{We only needed to write $x_0$ instead of $x$ to derive the formula for $f^*$. For simplicity we omit the index now.} Using $x=\frac{p}{2}$ we obtain
$$f^*(p) = px - f(x) = p\cdot \frac{p}{2} - f(p) = \frac{p^2}{4}.$$
It is not surprising that this coincides with the envelope $e(p)$ considered above, because there we started with the equation $x^2-px+q=0$ and regarded $q = px - x^2$ as a function of $p$. That is precisely what Legendre transformation for $f(x)=x^2$ does. The concept of Legendre transformation is important in physics; for further explanation and a more detailed introduction see \cite{Zia}.

\end{document}